%% file: convex_prog_robust_sparse_mean.tex
\documentclass[11pt, a4paper]{article}

\input{packages_notes}

\input{commandes_guillaume}

\DeclareMathOperator*{\Med}{Med}
\DeclareMathOperator*{\Tr}{Tr}

\usepackage{braket}
\usepackage{authblk}

\begin{document}
\title{Optimal robust mean and location estimation via convex programs with respect to any pseudo-norms}
\author[1]{Jules Depersin and Guillaume Lecu{\'e}  \\ email: \href{mailto:jules.depersin@ensae.fr}{jules.depersin@ensae.fr}, email: \href{mailto:lecueguillaume@gmail.com}{guillaume.lecue@ensae.fr} \\ CREST, ENSAE, IPParis. 5, avenue Henry Le Chatelier, 91120 Palaiseau, France.}


\maketitle

\begin{abstract}
We consider the problem of robust mean and location estimation w.r.t. any pseudo-norm of the form $x\in\bR^d\to\norm{x}_S = \sup_{v\in S}\inr{v,x}$ where $S$ is any symmetric subset of $\bR^d$. We  show  that the deviation-optimal minimax subgaussian rate for confidence $1-\delta$ is
\begin{equation*}
 \max\left(\frac{\ell^*(\Sigma^{1/2}S)}{\sqrt{N}},  \sup_{v\in S}\norm{\Sigma^{1/2}v}_2\sqrt{\frac{\log(1/\delta)}{N}}\right)
\end{equation*}where $\ell^*(\Sigma^{1/2}S)$ is the Gaussian mean width of $\Sigma^{1/2}S$ and $\Sigma$ the covariance of the data (in the benchmark i.i.d. Gaussian case). This improves the entropic minimax lower bound from \cite{MR4026610} and closes the gap characterized by Sudakov's inequality between the entropy and the Gaussian mean width for this problem. This shows that the right statistical complexity measure for the mean estimation problem is the Gaussian mean width. We also show that this rate can be achieved by a solution to a convex optimization problem in the adversarial and $L_2$ heavy-tailed setup by considering minimum of some Fenchel-Legendre transforms constructed using the Median-of-means principle. We finally show that this rate may also be achieved in situations where there is not even a first moment but a location parameter exists.
\end{abstract}


\section{Introduction} 
\label{sec:introduction}
We consider the problem of robust (to adversarial corruption and heavy-tailed data)  multivariate mean and location estimation with respect to any pseudo-norm $\nu\in\bR^d\to\norm{\nu}_S = \sup_{\mu\in S}\inr{\mu, \nu}$ where $S$ is any symmetric subset of $\bR^d$ (i.e. if $x\in S$ then $-x\in S$). This problem has been extensively studied during the last decade for $S=B_2^d$ the unit euclidean ball \cite{minsker2015geometric,catoni2018dimension,MR3491112,catoni2017dimension,devroye2016,lugosi2019sub,MR3052407,holland2019distribution,lerasle2019monk,MR4102693,Bartlett19,Jules_Guillaume_1,lei2020fast}. Only little is known for general symmetric sets $S$ and we will mainly refer to \cite{MR4026610} where this problem has been handled for $S$ which is the unit dual ball $B^\circ$ of a norm $\norm{\cdot}$ (so that $\norm{\cdot}_S = \norm{\cdot}$).

In \cite{MR4026610}, the authors introduced the problem of robust to heavy-tailed data estimation of a mean vector w.r.t. any norm. The problem can be stated as follow: given $N$ i.i.d. random vectors $ X_1, \ldots,  X_N$ in $\bR^d$ with mean $\mu^*$ and covariance matrix $\Sigma$, a norm $\norm{\cdot}$ on $\bR^d$ and a confidence parameter $\delta\in(0,1)$  find an estimator $\tilde \mu_N(\delta)$ and the best possible accuracy $r^*(N,\delta)$ such that with probability at least $1-\delta, \norm{\tilde \mu_N(\delta) - \mu^*}\leq r^*(N,\delta)$. In \cite{MR4026610}, the authors use the median-of-means principle \cite{MR702836,MR855970,MR1688610} to construct an estimator  satisfying the following result.

\begin{Theorem}\label{theo:lugosi_mendelson_1}[Theorem~2 in \cite{MR4026610}]
There exist an absolute constant $c$ such that the following holds. Given a norm $\norm{\cdot}$ on $\bR^d$ and a confidence $\delta\in(0,1)$, one can construct $\tilde \mu_N(\delta)$ such that with probability at least $1-\delta$
\begin{equation*}
\norm{\tilde \mu_N(\delta) - \mu^*}\leq \frac{c}{\sqrt{N}}\left(\bE \norm{\frac{1}{\sqrt{N}}\sum_{i=1}^N\eps_i (X_i - \mu^*)} + \bE \norm{\Sigma^{1/2}G} + \sup_{v\in B^\circ}\norm{\Sigma^{1/2}v}_2\sqrt{\log(1/\delta)}\right)
\end{equation*}where $B^\circ$ is the unit dual ball associated with $\norm{\cdot}$, $(\eps_i)$ are i.i.d. Rademacher variables independent of the $X_i$'s and $G\sim \cN(0,I_d)$.
\end{Theorem}
The construction of $\tilde \mu_N(\delta)$ is pretty involved and it seems hard to design an algorithm out of this procedure. In particular, $\tilde \mu_N(\delta)$ has not been proved to be solution to a convex optimization problem. Theorem~\ref{theo:lugosi_mendelson_1}'s main interest is thus from a theoretical point of view, while robust multivariate mean estimation can also be interesting from a practical point of view \cite{diakonikolas2017being}. 

The rate obtained in  Theorem~\ref{theo:lugosi_mendelson_1} can be decomposed into two terms: a deviation term $$\sup_{v\in B^\circ}\norm{\Sigma^{1/2}v}_2\sqrt{\log(1/\delta)}$$ where $\sup_{v\in B^\circ}\norm{\Sigma^{1/2}v}_2$ is a weak variance term and a complexity term which is the sum of a Rademacher complexity $\bE \norm{N^{-1/2}\sum_{i=1}^N\eps_i (X_i - \mu^*)}$ and a Gaussian mean width $\bE \norm{\Sigma^{1/2}G}$. The intuition behind this rate is explained in \cite{MR4026610}, in particular, in Question~1. We will however show that this rate is not the right one and that the Gaussian mean width term is actually not necessary. Moreover, we will show that the improved rate can be achieved by an estimator solution to a convex optimization problem in Section~\ref{sec:convex_programs} and that this holds even in the adversarial corruption model (see Assumption~\ref{assum:first} in Section~\ref{sec:convex_programs} below for a formal definition) and even in some situations where there is not even a first moment; in that case, $\mu^*$ is a \textit{location} parameter and $\Sigma$ a \textit{scatter} parameter.

The optimality of the rate in Theorem~\ref{theo:lugosi_mendelson_1} has been raised in \cite{MR4026610}. The classical approach to answer this type of  question is to consider the Gaussian case that is when the data $X_i,i\in[N]$ are i.i.d. $\cN(\mu^*, \Sigma)$. This is also the strategy used in \cite{MR4026610} to obtain the following deviation-minimax lower bound result\footnote{the result from \cite{MR4026610} is proved  for $\Sigma = I_d$, it is however straightforward to extend it to the general case.}. 

\begin{Theorem}\label{theo:lugosi_mendelson_lower_bound}[Theorem~3 and first paragraph in p.962 in \cite{MR4026610}]
There exists an absolute constant $c>0$ such that the following holds.
If $\hat \mu:\bR^{Nd}\to \bR^d$ is an estimator such that for all $\mu^*\in\bR^d$ and all $\delta\in(0,1/4)$,
\begin{equation*}
\bP_{\mu^*}^N\left[\norm{\hat \mu - \mu^*}\leq r^*\right]\geq 1-\delta
\end{equation*}where $\bP_{\mu^*}^N$ is the probability distribution of $(X_i)_{i\in[N]}$ when the $ X_i$ are i.i.d. $\cN(\mu^*, \Sigma)$ then
\begin{equation*}
r^*\geq \frac{c}{\sqrt{N}}\left( \sup_{\eta>0} \eta \sqrt{\log N(\Sigma^{1/2} B^\circ, \eta B_2^d)} + \sup_{v\in B^\circ}\norm{\Sigma^{1/2}v}_2\sqrt{\log(1/\delta)}\right)
\end{equation*}where $N(\Sigma^{1/2} B^\circ, \eta B_2^d)$ is the minimal number of translated of $\eta B_2^d$ needed to cover $\Sigma^{1/2} B^\circ$.
\end{Theorem}

The term $\sup_{v\in S}\norm{\Sigma^{1/2}v}_2\sqrt{\log(1/\delta)}$ in the lower bound from Theorem~\ref{theo:lugosi_mendelson_lower_bound} is obtained in \cite{MR4026610} from Proposition~6.1 in \cite{MR3052407} which is a deviation-minimax lower bound result holding in the one dimensional case which relies on the fact that the empirical mean is a sufficient statistics in the Gaussian shift theorem\footnote{The argument used in \cite{MR4026610} goes from the one dimensional case studied in \cite{MR3052407} to the $d$-dimensional case. It is given in a none formal way and may require some extra argument to hold. Indeed the estimator $x^*(\hat \Psi_N)$ in \cite{MR4026610} is constructed using the $d$-dimensional data $X_1, \ldots, X_N$ and not one-dimensional data such as $x^*(X_1), \ldots, x^*(X_N)$. However, the result from \cite{MR3052407} holds for estimators of a one dimensional mean using one-dimensional data and not $d$-dimensional ones. Nevertheless, Olivier Catoni showed us how to adapt the proof of Proposition~6.1 in \cite{MR3052407} by using the sufficiency of the empirical mean in the Gaussian shift model in  $\bR^d$ to get this deviation dependent lower bound term.}.  

The complexity term $\sup_{\eta>0} \eta \sqrt{\log N(\Sigma^{1/2} B^\circ, \eta  B_2^d)}$ obtained in Theorem~\ref{theo:lugosi_mendelson_lower_bound} follows from the duality theorem of metric entropy from \cite{MR2113023} and a volumetric argument in the Gauss space similar to the one used to prove dual Sudakov's inequality in p.82-83 in \cite{MR2814399} which has also been used to obtain minimax lower bounds based on the entropy in \cite{lecue2013learning} and \cite{MR3718154}.

In general, there is a gap between the upper bound from Theorem~\ref{theo:lugosi_mendelson_1} and the lower bound from Theorem~\ref{theo:lugosi_mendelson_lower_bound} even in the Gaussian case. This gap is characterized by Sudakov's inequality (see Theorem~3.18 in \cite{MR2814399} or Theorem~5.6 in \cite{MR1036275}):
\begin{equation}\label{eq:sudakov_bound}
\sup_{\eta>0} \eta \sqrt{\log N(\Sigma^{1/2} B^\circ, \eta B_2^d)} \leq c \bE \norm{\Sigma^{1/2}G}
\end{equation}where $G\sim\cN(0,I_d)$. Indeed, in the Gaussian case the complexity term of the rate obtained in Theorem~\ref{theo:lugosi_mendelson_1} is the Gaussian mean width, that is the right-hand term from \eqref{eq:sudakov_bound} whereas the complexity term from Theorem~\ref{theo:lugosi_mendelson_lower_bound} is the entropy, that is the left-hand term in \eqref{eq:sudakov_bound}.

As mentioned in Remark~3 from \cite{MR4026610}, when Sudakov's inequality \eqref{eq:sudakov_bound} is sharp then upper and lower bounds from Theorem~\ref{theo:lugosi_mendelson_1} and \ref{theo:lugosi_mendelson_lower_bound} match in the Gaussian case (in that case the Rademacher complexity is equal to the Gaussian mean width in Theorem~\ref{theo:lugosi_mendelson_1}). Sharpness in Sudakov's inequality is however not a typical situation. In particular, for ellipsoids, Sudakov's bound \eqref{eq:sudakov_bound} is not sharp in general and therefore the lower bound from Theorem~\ref{theo:lugosi_mendelson_lower_bound} fails to recover the classical  subgaussian rate for the standard Euclidean norm case (that is for $S=B_2^d$) which is given in \cite{lugosi2019sub} by
\begin{equation}\label{eq:subgaussian_rate}
\sqrt{\frac{\Tr{(\Sigma)}}{N}} + \sqrt{\frac{\norm{\Sigma}_{op}\log(1/\delta)}{N}}.
\end{equation} Indeed, when $\norm{\cdot}$ is the $\ell_2^d$ Euclidean norm then $\bE \norm{\Sigma^{1/2}G}=\bE \norm{\Sigma^{1/2}G}_2\sim \sqrt{\Tr(\Sigma)}$ (see, for instance, Proposition~2.5.1 in \cite{MR3184689}). Whereas, for the entropy of $\Sigma^{1/2} B^\circ=\Sigma^{1/2} B_2^d$ w.r.t. $\eta B_2^d$, it follows from equation~(5.45) in \cite{MR1036275} that 
\begin{equation}\label{eq:entropy_ellipse_pisier}  \sup_{\eta>0} \eta \sqrt{\log_2 N(\Sigma^{1/2}B_2^d, \eta B_2^d)} = \sup_{n\geq1} e_{n+1}(\Sigma^{1/2})\sqrt{n+1}\sim \sup_{n\geq1,k\in[d]} \frac{\sqrt{n}}{2^{n/k}}\left|\prod_{j=1}^k \sqrt{\lambda_j}\right|^{1/k} \sim \sqrt{\sup_{k\in[d]}k\left|\prod_{j=1}^k \lambda_j\right|^{1/k}}
\end{equation}where $(e_{n+1}(\Sigma^{1/2}))_n$ are the entropy numbers of $\Sigma^{1/2}:\ell_2^d\to \ell_2^d$ (see page~62 in \cite{MR1036275} for a definition) and $\lambda_1\geq\ldots\geq \lambda_d$ are the singular values of $\Sigma$. In particular, when $\lambda_j = 1/j$, the entropy bound \eqref{eq:entropy_ellipse_pisier} is of the order of a constant whereas the Gaussian mean width is of the order of $\sqrt{\log d}$. We will fill this gap in Section~\ref{sec:benchmark_results_subgaussian_rates} by showing a lower bound where the entropy is replaced by the (larger) Gaussian mean width. We will therefore obtain matching upper and lower bounds revealing that Gaussian mean width is the right way to measure the statistical complexity for the mean estimation problem w.r.t. any $\norm{\cdot}_S$.

The paper is organized as follows. In the next section, we obtain the deviation-minimax optimal rate in the i.i.d. Gaussian case. In Section~\ref{sec:convex_programs} we show that the rate from Theorem~\ref{theo:lugosi_mendelson_1} can be improved and that it can be achieved by a solution to a convex program in the adversarial contamination model and in under weak or no moment assumptions. All the proofs have been gathered in Section~\ref{sec:proofs}.


\section{Deviation minimax rates in the Gaussian case: benchmark subgaussian rates for the mean estimation w.r.t. $\norm{\cdot}_S$} 
\label{sec:benchmark_results_subgaussian_rates}
In this section, we obtain the optimal deviation-minimax rates of estimation of a mean vector $\mu^*$ when we are given $N$ i.i.d. $X_1, \ldots, X_N$ distributed like $\cN(\mu^*, \Sigma)$ when $\Sigma\succeq0$ is some unknown covariance matrix. In the following, $\bP_{\mu^*}^N$ denotes the probability distribution of $(X_1, \ldots, X_N)$; it is a Gaussian measure on $\bR^{Nd}$ with mean $((\mu^*)^\top, \ldots, (\mu^*)^\top)$ and a block $(Nd)\times (Nd)$ covariance matrix with $d\times d$ diagonal blocks given  by $\Sigma$ repeated $N$ times and $0$ outside of these blocks.

Unlike classical minimax results holding in expectation or with constant probability (see Chapter~2 in \cite{MR2724359}) we want, in this section, the deviation parameter $\delta$ to appear explicitly in the minimax lower bound. Moreover, this dependency of the convergence rate with respect to $\delta$ should be of the right order given by the subgaussian $\sqrt{\log(1/\delta)}$ rate and not other polynomial dependency such as $\sqrt{1/\delta}$ as one gets for the empirical mean for $L_2$ variables (see Proposition~6.2 in \cite{MR3052407}). This subtle behavior of the rate in terms of $\delta$ cannot be seen in expectation or constant deviation minimax lower bounds. In particular,  this makes such results (like Theorem~\ref{theo:minimax_gauss_S} or \ref{theo:minimax_gauss} below) unachievable via classical information theoretic arguments as in Chapter~2 in \cite{MR2724359}. 

Fortunately, in \cite{lecue2013learning}, a minimax lower bound has been proved thanks to the Gaussian shift theorem which makes the deviation parameter $\delta$ appearing explicitly in the minimax lower bound. We use the same strategy here to prove our main result Theorem~\ref{theo:minimax_gauss_S} below and its corollary Theorem~\ref{theo:minimax_gauss} in the classical Euclidean $S=B_2^d$ case.

We consider the general problem of estimating $\mu^*$ w.r.t. $\norm{\cdot}_S$. Let $S\subset \bR^d$ be a symmetric set. We first obtain an upper bound result revealing the subgaussian rate. We use the empirical mean $\bar X_N=N^{-1}\sum_i X_i$ as an estimator of $\mu^*$. Using Borell TIS's inequality (Theorem~7.1 in \cite{Led01} or pages 56-57 in \cite{MR3184689}) we get: for all $0<\delta <1$, with probability at least $1-\delta$,
\begin{equation*}
\norm{\bar X_N-\mu}_S = \sup_{v\in S}\inr{v, \bar X_N-\mu}\leq \bE \sup_{v\in S}\inr{v, \bar X_N-\mu} + \sigma_S \sqrt{2\log(1/\delta)}
\end{equation*}where $\sigma_S = \sup_{v\in S} \sqrt{\bE \inr{v, \bar X_N-\mu}^2}$ is called the weak variance. It follows that with probability at least $1-\delta$,
\begin{equation}\label{eq:subgauss_rate_S}
\norm{\bar X_N-\mu}_S\leq \frac{\ell^*(\Sigma^{1/2}S)}{\sqrt{N}} + \frac{\sup_{v\in S}\norm{\Sigma^{1/2}v}_2\sqrt{\log(1/\delta)} }{\sqrt{N}}
\end{equation}where $\ell^*(\Sigma^{1/2}S) = \sup\big(\inr{G, x}: x\in \Sigma^{1/2}S\big) = \bE \norm{\Sigma^{1/2}G}_S$, for $G\sim \cN(0, I_d)$, is the Gaussian mean width of the set $\Sigma^{1/2}S$. In particular, in the case where $S=B_2^d$, we recover the subgaussian rate  \eqref{eq:subgaussian_rate} in \eqref{eq:subgauss_rate_S}. Our aim is now to show that the rate in \eqref{eq:subgauss_rate_S} is deviation-minimax optimal. This is what is obtained in the next result.

\begin{Theorem}\label{theo:minimax_gauss_S} Let $S$ be a symmetric subset of $\bR^d$ such that ${\rm span}(S) = \bR^d$. If $\hat \mu:\bR^{Nd}\to \bR^d$ is an estimator such that for all $\mu^*\in\bR^d$ and all $\delta\in(0,1/4]$,
\begin{equation*}
\bP_{\mu^*}^N\left[\norm{\hat \mu - \mu^*}_S\leq r^*\right]\geq 1-\delta
\end{equation*}then
\begin{equation*}
r^*\geq \max\left(\frac{1}{24}\sqrt{\frac{\log 2}{\log(5/4)}}\frac{\ell^*(\Sigma^{1/2}S)}{\sqrt{N}},  \frac{\sup_{v\in S}\norm{\Sigma^{1/2}v}_2}{12}\sqrt{\frac{\log(1/\delta)}{\sqrt{N}}}\right).
\end{equation*}
\end{Theorem}

It follows from the upper bound \eqref{eq:subgauss_rate_S} and the deviation-minimax lower bound from Theorem~\ref{theo:minimax_gauss_S} that it is now possible to know exactly (up to absolute constants) the subgaussian rate for the problem of mean estimation in $\bR^d$ w.r.t. $\norm{\cdot}_S$, it is given by
\begin{equation}\label{eq:subgaussian_rate_S}
 \max\left(\frac{\ell^*(\Sigma^{1/2}S)}{\sqrt{N}},  \frac{\sup_{v\in S}\norm{\Sigma^{1/2}v}_2\sqrt{\log(1/\delta)} }{\sqrt{N}}\right).
 \end{equation} We may identify the two complexity and deviation terms in this rate. In particular, the complexity term is measured here via the Gaussian mean width of the set $\Sigma^{1/2}S$ and not its entropy as it was previously known following Theorem~\ref{theo:lugosi_mendelson_lower_bound}. Theorem~\ref{theo:minimax_gauss_S} together with \eqref{eq:subgauss_rate_S} show that the right way to measure the statistical complexity in the problem of mean estimation in $\bR^d$ w.r.t. to any $\norm{\cdot}_S$ is via the Gaussian mean width. This differs from other statistical problems such as the regression model with random design where the entropy has been proved to be the right statistical complexity in several examples \cite{MR3718154,lecue2013learning}. Following the later results in the regression model, Theorem~\ref{theo:minimax_gauss_S} is a bit unexpected because one may though that by taking an ERM over an epsilon net of $\bR^d$ for the right choice of $\eps$ one could obtain a better rate than the one driven by the Gaussian mean width in \eqref{eq:subgaussian_rate_S}; indeed, for this type of procedure, one may expect a rate depending on the (smaller) entropy instead of the (larger) Gaussian mean width. Theorem~\ref{theo:minimax_gauss_S} shows that this is not the case: even discretized ERM cannot achieve a better rate than the one driven by the Gaussian mean width in the mean estimation problem.

An important consequence of Theorem~\ref{theo:minimax_gauss_S} is obtained when $S=B_2^d$ that is for the problem of multivariate mean estimation w.r.t. the $\ell_2^d$-norm which is the problem that has been extensively considered during the last decade.  In the following result, we recover the well-known subgaussian rate \eqref{eq:subgaussian_rate} showing that all the upper bound results where this rate has been proved to be achieved are actually deviation-minimax optimal and therefore could not have been improved uniformly over all $\mu^*\in\bR^d$.

\begin{Theorem}\label{theo:minimax_gauss}
If $\hat \mu:\bR^{Nd}\to \bR^d$ is an estimator such that $\bP_{\mu^*}^N\left[\norm{\hat \mu - \mu^*}_2\leq r^*\right]\geq 1-\delta$ for all $\mu^*\in\bR^d$ and all $\delta\in(0,1/4]$, then
\begin{equation*}
r^*\geq \max\left(\frac{1}{24}\sqrt{\frac{\log 2}{2\log(5/4)}}\sqrt{\frac{{\rm Tr}(\Sigma)}{N}},  \frac{1}{12}\sqrt{\frac{\norm{\Sigma}_{op}\log(1/\delta)}{N}}\right).
\end{equation*}
\end{Theorem}

Given that the empirical mean $\bar X_N$ is such that for all $\mu\in\bR^d$ with $\bP_\mu^N$-probability at least $1-\delta$,
\begin{equation*}
\norm{\bar X_N - \mu}_2\leq \sqrt{\frac{\Tr{(\Sigma)}}{N}} + \sqrt{\frac{2\norm{\Sigma}_{op}\log(1/\delta)}{N}}
\end{equation*}we conclude from Theorem~\ref{theo:minimax_gauss} that the sub-gaussian rate \eqref{eq:subgaussian_rate} is the deviation-minimax rate of convergence for the multivariate mean estimation problem w.r.t. $\ell_2^d$ and that it is achieved by the empirical mean. In particular, there are no statistical procedure that can do better than the empirical mean uniformly over all mean vectors $\mu^*\in\bR^d$ up to constant, this includes in particular all discretized versions of $\bar X_N$.


\section{Convex programs} 
\label{sec:convex_programs}
In this section, we introduce statistical procedures which are solutions to convex programs and which can achieve the rate from Theorem~\ref{theo:lugosi_mendelson_1} without the unnecessary Gaussian mean width term $\bE\norm{\Sigma^{1/2}G}$. We also show that these procedures handle adversarial corruption and may still perform optimally in some situations where there is not even a first moment.
\subsection{Construction of the Fenchel-Legendre minimum estimators.} 
\label{par:construction_of_the_fenchel_legendre_minimum_estimators_}

\begin{Definition}
Let $S$ be a subset of $\bR^d$ and $f:\bR^d\to \bR$. The Fenchel-Legendre transform of $f$ on $S$ is the function $f^*_S$ defined for all $\mu\in\bR^d$ by $f_S^*(\mu) = \sup_{v\in S} \left(\inr{\mu, v} - f(v)\right)$.
\end{Definition}
For our purpose, the main property of a Fenchel-Legendre transform we will use is that it is a convex function as it is the maximal function of the family $(\mu\in\bR^d\to \inr{\mu, v} - f(v): v\in S)$ of linear functions. 

We are now defining two examples of functions such that by taking the minimum of their Fenchel-Legendre transform over $S$ will lead to optimal estimators of $\mu^*$ w.r.t. $\norm{\cdot}_S$. The construction of these two functions are based on the median-of-means principle: the dataset $\{X_1, \ldots, X_N\}$ is split into $K$ equal size blocks of data indexed by $(B_k)_{k}$ forming an equipartition of $[N]$. On each block, an empirical mean is constructed $\bar X_k = |B_k|^{-1}\sum_{i\in B_k}X_i$.  The two functions we are considering are using the $K$ bucketed means $(\bar X_k)_k$ and are defined, for all $v\in \bR^d$, by
\begin{equation}\label{eq:two_main_functions}
f(v) = \frac{1}{|I_K|}\sum_{k\in I_K} \inr{\bar X_k, v}_{(k)}^*   \mbox{ and } g(v) =\Med(\inr{\bar X_k,v}) = \inr{\bar X_k,v}_{\left(\frac{K+1}{2}\right)}^*
\end{equation}where if $a_k=\inr{\bar X_k, v}, k\in[K]$ then $\inr{\bar X_k, v}_{(k)}^*,k\in[K]$ are the rearrangement of $(a_k)_k$ such that $a_{(1)}^*\leq \ldots \leq a_{(K)}^*$ (this is the rearrangement of the values $a_k$'s themselves and not of their absolute values)  and 
\begin{equation*}
I_K=\left[\frac{K+1}{4}, \frac{3(K+1)}{4}\right] = \left\{\frac{K+1}{2}\pm k:k=0,1,\cdots, \frac{K+1}{4}\right\}
\end{equation*}is the inter-quartiles interval -- w.l.o.g. we assume that $K+1$ can be divided by $4$. In other words, $f(v)$ is the average sum over all inter-quartile values of the vector $(\inr{\bar X_k, v})_{k\in[K]}$ and $g(v)$ is the median of this vector.  Note that both functions $f$ and $g$ are homogeneous i.e. $f(\theta v)=\theta f(v)$ and $g(\theta v)=\theta g(v)$ for every $v\in\bR^d$ and $\theta\in\bR$ and in particular they are odd functions; two facts we will use later. 

We are now considering the Fenchel-Legendre transform of the functions $f$ and $g$ over a symmetric set $S$:
\begin{equation}\label{eq:legendre_fenchel_transfrom}
f^*_S:\mu\in\bR^d \to\sup_{v\in S} \left(\inr{\mu, v} - f(v)\right) \mbox{ and } g^*_S:\mu\in\bR^d \to\sup_{v\in S} \left(\inr{\mu, v} - g(v)\right). 
\end{equation} As mentioned previously the two functions $f^*_S$ and $g_S^*$ are convex functions. We are now using them to define convex programs whose solutions will be proved to be robust and subgaussian estimators of the mean / location vector $\mu^*$ w.r.t. $\norm{\cdot}_S$:
\begin{equation}\label{eq:estimators_basics}
\hat \mu_S^f \in\argmin_{\mu\in\bR^d} f^*_S(\mu) \mbox{ and } \hat \mu_S^g \in\argmin_{\mu\in\bR^d} g^*_S(\mu). 
\end{equation}

For some special choices of $S$, the Fenchel-Legendre minimization estimator $\hat \mu_S^g$ coincides with some classical procedures. This is for instance the case when $S=B_1^d$ (the unit ball of the $\ell_1^d$-norm) or $S=B_2^d$. Indeed, when $S=B_1^d$, $\hat \mu_S^g$ is the coordinate-wise Median of Means:
\begin{equation}\label{eq:coordinate_wise_MOM}
\hat\mu_S^g=\argmin_{\mu=(\mu_j)\in\bR^d}\max_{j\in[d]}\left|\mu_j - \Med\left(\inr{\bar X_k, e_j}\right)\right| = \left(\Med\left(\inr{\bar X_k, e_j}\right):j\in[d]\right)
\end{equation}where $(e_j)_{j=1}^d$ is the canonical basis of $\bR^d$, because $\norm{\cdot}_S = \norm{\cdot}_{{\rm conv}(S)}$ where ${\rm conv}(S)$ is the convex hull of $S$ and so one may just take  $S=\{\pm e_j:j\in[d]\}$. It is therefore possible to derive deviation-minimax optimal bounds for the coordinate-wise Median of Means w.r.t. the $\ell_\infty^d$-norm from general upper bounds on $\hat \mu^g_S$ since in that case $\norm{\cdot}_S=\norm{\cdot}_\infty$.

In the case $S=B_2^d$ (that is for the mean/location estimation problem w.r.t. $\ell_2^d$), the Fenchel-Legendre minimum  estimator $\hat \mu_S^g$ is a minmax MOM estimator \cite{lecue2020}. This connection allows to write $\hat \mu_S^g$ (as well as $\hat \mu_S^f$)  as a non-constraint estimator, it also shows that this minmax MOM estimator is actually solution to a convex optimization problem and how minmax MOM estimator can be generalized to other estimation risks.

Minmax MOM estimators have been introduced as a systematic way to construct robust and subgaussian estimators in  \cite{lecue2020}. They have been proved to be deviation-minimax optimal for the mean estimation problem in \cite{lerasle2019monk} w.r.t. $\norm{\cdot}_2$. Their definition only requires to consider a loss function; here we take for all $\mu\in\bR^d$, $\ell_\mu:x\in\bR^d\to \norm{x-\mu}_2^2$ and the minmax MOM estimator is then defined as 
\begin{equation}\label{eq:minmax_MOM}
\tilde\mu\in\argmin_{\mu\in\bR^d}\sup_{\nu\in\bR^d} \Med\left(P_{B_k}(\ell_\mu - \ell_\nu):k\in[K]\right)
\end{equation}where $P_{B_k}$ is the empirical measure on the data in block $B_k$.  The minmax MOM estimator $\tilde\mu$ was proved to achieve the subgaussian rate in \eqref{eq:subgaussian_rate} with confidence $1-\delta$ when the number of blocks is $K\sim \log(1/\delta)$ and $K\gtrsim |\cO|$ in \cite{lerasle2019monk}. 

Even though the minmax formulation of $\tilde\mu$ suggests a robust version of a descent/ascent gradient method over the median block (see \cite{lecue2020,lerasle2019monk} for more details), no proof of convergence of this algorithm is known so far. Moreover, the main drawback of the minmax MOM estimator seems to be that it is solution of a non-convex optimization problem and may therefore be likely to be rather difficult to compute in practice. In the next result, we show that this is not the case since the minmax MOM estimator \eqref{eq:minmax_MOM} is in fact equal to $\hat \mu_S^g$ for $S=B_2^{d}$ and it is therefore solution to a convex optimization problem.

\begin{Proposition}\label{prop:minmax_mom_convex}
The minmax MOM estimator $\tilde \mu$ defined in \eqref{eq:minmax_MOM} satisfies  $\tilde \mu \in\argmin_{\mu\in\bR^d} g^*_{B_2^{d}}(\mu)$. The minmax MOM estimator is therefore solution to a convex optimization problem. 
\end{Proposition}

\beginproof We show that $\tilde \mu \in \argmin_{\mu\in\bR^d} \sup_{\norm{v}_2=1}\Med(\inr{\bar X_k - \mu, v})$. We consider the quadratic/multiplier decomposition of the difference of loss functions: for all $\mu,\nu\in\bR^d$ and $x\in\bR^d$, we have $(\ell_\mu - \ell_\nu)(x) = \norm{x-\mu}_2^2- \norm{x-\nu}_2^2 = -2\inr{x-\mu, \mu-\nu} - \norm{\mu-\nu}_2^2$. Hence, for all $\mu\in\bR^d$, we have
\begin{align*}
&\sup_{\nu\in\bR^d}\Med\left(P_{B_k}(\ell_\mu - \ell_\nu)\right) = \sup_{\nu\in\bR^d}\left(-2\Med(\inr{\bar X_k - \mu, \mu-\nu}) - \norm{\mu-\nu}_2^2\right)\\
& = \sup_{\norm{v}_2=1}\sup_{\theta\geq0} \left(2\theta\Med(\inr{\bar X_k - \mu, v}) - \theta^2\right) = \sup_{\norm{v}_2=1}\left(\Med(\inr{\bar X_k - \mu, v})\right)^2 = \left(\sup_{\norm{v}_2=1}\Med(\inr{\bar X_k - \mu, v})\right)^2.
\end{align*}We conclude since
\begin{equation*}
\argmin_{\mu\in\bR^d}\left(\sup_{\norm{v}_2=1}\Med(\inr{\bar X_k - \mu, v})\right)^2 = \argmin_{\mu\in\bR^d} \sup_{\norm{v}_2=1}\Med\left(\inr{\bar X_k - \mu, v}\right).
\end{equation*}
\endproof

It follows from Proposition~\ref{prop:minmax_mom_convex} that the minmax MOM estimator $\tilde\mu$ is solution to  a convex optimization problem. This fact is far from being obvious given the definition of $\tilde\mu$ in \eqref{eq:minmax_MOM}.


Proposition~\ref{prop:minmax_mom_convex} suggests a new formulation for $\hat \mu_S^g$ and $\hat\mu_S^f$. It is indeed possible to write these estimators as regularized estimators instead of their original constraint formulation (note that the Fenchel-Legendre transforms in \eqref{eq:legendre_fenchel_transfrom} are suprema over $S$ and are therefore constraint optimization problems). We now show that we may write them as suprema over all $\bR^d$ if we add an ad hoc regularization function.

Let us introduce the two following functions which may be seen as regularized versions of the two $f$ and $g$ functions from \eqref{eq:two_main_functions}: for all $\nu\in\bR^d$,
\begin{equation}\label{eq:fct_F}
F_S(\nu) = f(\nu) + \frac{\norm{\nu}_S^2}{4} \mbox{ and } G_S(v) = g(\nu)+ \frac{\norm{\nu}_S^2}{4}.
  \end{equation} We also consider their Fenchel-Legendre transforms over the entire set $\bR^d$: for all $\mu\in\bR^d$,
 \begin{equation*}
 F_S^*(\mu) = \sup_{\nu\in\bR^d}\left(\inr{\mu, \nu} - F_S(\nu)\right) \mbox{ and } G_S^*(\mu) = \sup_{\nu\in\bR^d}\left(\inr{\mu, \nu} - G_S(\nu)\right).
 \end{equation*}

 The next result shows that the later two Fenchel-Legendre transforms can be used to define the two estimators $\hat\mu_S^f$ and $\hat\mu_S^g$. The proof of Proposition~\ref{prop:F_and_G} is similar to the one of Proposition~\ref{prop:minmax_mom_convex} where the $\ell_2$-norm is replaced by $\norm{\cdot}_S$ and is therefore omitted.

 \begin{Proposition}\label{prop:F_and_G}
 Let $S$  be a symmetric subset of $\bR^d$ such that ${\rm span}(S)=\bR^d$. We have $\hat \mu_S^f \in\argmin_{\mu\in\bR^d} F^*_S(\mu)$  and   $\hat \mu_S^g \in\argmin_{\mu\in\bR^d} G^*_S(\mu)$.
 \end{Proposition}

 As a consequence of Proposition~\ref{prop:F_and_G}, one can write the two estimators $\hat \mu_S^f$ and $\hat \mu_S^g$ as solutions to unconstrained minmax optimization problems like the minmax MOM estimator \eqref{eq:minmax_MOM} and in particular, one may design an alternating ascent/descent sub-gradient algorithm similar to the one from \cite{lecue2020} -- we expect the one associated with $\hat \mu_S^f$ which uses half of the dataset at each iteration to be more efficient than the one associated with $\hat \mu_S^g$ which uses only the $N/K$ data in the median block at each iteration. That is the reason why we provide in Figure~\ref{algo:ascent_descent} this algorithm only for
 \begin{equation*}
  \hat \mu_S^f\in\argmin_{\mu\in\bR^d} \sup_{\nu\in\bR^d} \left(\inr{\mu, \nu} - \frac{1}{|I_K|}\sum_{k\in I_K} \inr{\bar X_k, v}_{(k)}^* - \frac{\norm{\nu}_S^2}{4} \right).
  \end{equation*} 

 \vspace{0.7cm}
 \begin{algorithm}[H]\label{algo:ascent_descent}
\SetKwInOut{Input}{input}\SetKwInOut{Output}{output}\SetKw{Or}{or}
\SetKw{Return}{Return}
\Input{the data $X_1, \ldots, X_N$, a number $K$ of blocks, two steps size sequences $(\eta_t)_t, (\theta_t)_t$  and $\eps>0$ a stopping parameter}
\Output{A robust estimator of the mean $\mu$} 
Construct an equipartition $B_1\sqcup \cdots \sqcup B_K=\{1,\cdots,N\}$ at random\\
Construct the $K$ empirical means $\bar{X}_k=(N/K)\sum_{i\in B_k}X_i, k\in[K]$\\
Compute $\tilde\mu^{(0)}$ the coordinate-wise median-of-means and put $\mu^{(0)} = \tilde\mu^{(0)}$ and $\nu^{(0)} = \tilde\mu^{(0)}$\\
\While{$\norm{\mu^{(t)} - \mu^{(t+1)}}_S\geq \eps$}{ 
\BlankLine
Construct an equipartition $B_1\sqcup \cdots \sqcup B_K=\{1,\cdots,N\}$ at random\\
Construct the $K$ empirical means $\bar{X}_k=(N/K)\sum_{i\in B_k}X_i, k\in[K]$\\
Find the inter-quartile block numbers $k_1, \ldots, k_{(K+1)/2}\in[K]$ such that
\begin{equation*}
  f(\nu^{(t)}) = \frac{1}{|I_K|}\sum_{j=1}^{(K+1)/2} \inr{\bar X_{k_j}, \nu^{(t)}}.
  \end{equation*}  
Construct $g^{(t)}$ a subgradient of $\norm{\cdot}_S$ at $\nu^{(t)}$ and the ascent direction
\begin{equation*}
\nabla_\nu^{(t+1)} = \mu^{(t)} - \frac{1}{|I_K|}\sum_{j=1}^{(K+1)/2} \bar X_{k_j} - \frac{\norm{\nu^{(t)}}_S g^{(t)}}{2}.
\end{equation*}
Update $\nu^{(t+1)}\leftarrow \nu^{(t)}+\eta_t \nabla_\nu^{(t+1)}$.\\
Make one descent step: $\mu^{(t+1)}\leftarrow \mu^{(t)}-\theta_t \nu^{(t+1)}$.\\}
\Return $\mu^{(t+1)}$
 \caption{An alternating ascent/descent algorithm for the robust mean estimation problem w.r.t. $\norm{\cdot}_S$ with randomly chosen  blocks of data at each step.}
\end{algorithm}
\vspace{0.7cm}


\subsection{The adversarial corruption model and two models for inlier.} 
\label{par:the_adversarial_corruption_with_}
In this section, we introduce the assumptions under which we will obtain some statistical upper bounds for the Fenchel-Legendre minimum estimators introduced above. We are considering two types of assumptions: one for the outliers which will be the adversarial corruption model and one for the inlier which will be either the existence of a second moment or a regularity assumption on a family of cdf around $0$. We start with the adversarial corruption model.

\begin{Assumption}\label{assum:first} There exists $N$ independent random vectors $(\tilde X_i)_{i=1}^N$ in $\bR^d$.  The $N$ random vectors $(\tilde X_i)_{i=1}^N$ are first given to an ''adversary'' who is allowed to modify up to $|\cO|$ of these vectors. This modification does not have to follow any rule. Then, the ''adversary'' gives the modified dataset $(X_i)_{i=1}^N$ to the statistician. Hence, the statistician receives an ''adversarially'' contaminated dataset of $N$ vectors in $\bR^d$ which can be partitioned into two groups: the modified data $(X_i)_{i\in\cO}$, which can be seen as outliers and the ''good data'' or inlier $(X_i)_{i\in\cI}$ such that $\forall i\in\cI, X_i=\tilde X_i$. Of course, the statistician does not know which data has been modified or not so that the partition $\cO\cup\cI=\{1, \ldots, N\}$ is unknown to the statistician.  
\end{Assumption}

In the adversarial contamination model from Assumption~\ref{assum:first}, the set  $\cO\subset [N]$ can depend arbitrarily on the initial data $(\tilde X_i)_{i=1}^N$;  the corrupted data $(X_i)_{i \in \cO}$ can have any arbitrary dependence structure; and the informative data $(X_i)_{i\in\cI}$ may also be correlated (for instance, it is, in general, the case  when the $|\cO|$ data $\tilde X_i$ with largest $\ell_2^d$-norm are modified by the adversary). The adversarial corruption model covers the Huber $\eps$-contamination model \cite{MR2488795} and also the $\cO\cup\cI$ framework from \cite{lecue2019learning,lecue2020,LMSL}. 

Assumption~\ref{assum:first} does not grant any property of the inlier data $(\tilde X_i)_{i\in[N]}$ except that they are independent. We will obtain a general result under only  Assumption~\ref{assum:first} in Section~\ref{sec:proofs}. However, to recover convergence rates similar to the one in  Theorem~\ref{theo:lugosi_mendelson_1} or the subgaussian rate in \eqref{eq:subgaussian_rate_S}, we will grant some assumptions on the inlier as well. We are now considering two assumptions on the inlier which are of different nature.

The two assumptions on the inlier we are now considering are related to a subtle property of the Median-of-Means (MOM) principle which somehow benefits from its two components: the empirical median and the empirical mean. Indeed,  MOM is en empirical median of empirical means and so if we refer to the classical asymptotic normality (a.n.) results of the empirical mean and the empirical median, the first one holds under the existence of a second moment and the second one holds under the assumption that the cdf is differentiable at the median with positive derivative at the median (see Corollary~21.5 in \cite{MR1652247}). We therefore recover these two types of assumptions when we work with estimators using the MOM principle. A nice feature of MOM based estimators is that their estimation results hold  under either one of the two conditions and do not require the two assumptions to hold simultaneously. We can therefore consider the two assumptions independently and get two estimation results for the Fenchel-Legendre minimum estimators introduced above (which are based on the MOM principle). We start with the moment assumption.

\begin{Assumption}\label{ass:moment}
The $N$ independent random vectors $(\tilde X_i)_{i=1}^N$ have mean $\mu^*$ and there exists a SDP matrix $\Sigma\in\bR^{d\times d}$ such that $\bE (\tilde X_i - \mu^*)(\tilde X_i - \mu^*)^\top\preceq \Sigma$. 
\end{Assumption}

Most of the statistical bounds obtained on MOM based estimators have focused on the heavy-tailed setup and have therefore consider  Assumption~\ref{ass:moment} as their main assumption. This is the 'empirical mean component' of the MOM principle which has been the most exploited so far. It is however also possible to use the 'empirical median component' of the MOM principle to get statistical bounds even in cases where a first moment does not even exist. In that case, $\mu^*$ is called a \textit{location parameter} and $\Sigma$ a \textit{scale parameter}. Also, a natural assumption is similar to the one used to get the a.n. of the empirical median, that is an assumption on the cdf at the median adapted to the multidimensional and non-asymptotic setup. We are now introducing such an assumption.

\begin{Assumption}\label{ass:function_H}
The inlier data $(\tilde X_i)_{i=1}^N$ are i.i.d.. There exists $\mu^*\in\bR^d$ and two absolute constants $c_0>0$ and $c_1>0$ such that the following holds: for all $v\in S$ and all $0<r\leq c_0$, $H_{N,K,v}(r)\leq 1/2-c_1 r$ where
\begin{equation}\label{eq:function_H}
H_{N,K,v}(r) = \bP\left[\frac{1}{\sqrt{N/K}}\sum_{i=1}^{N/K}\inr{\tilde X_i - \mu^*, v}> r\right].
\end{equation}
\end{Assumption}  

A typical example where Assumption~\ref{ass:function_H} holds is when $S=\cS_2^{d-1}$ (that is for the location estimation problem w.r.t. the Euclidean $\ell_2^d$ norm) and the $\tilde X_i$'s are rotational invariant that is when for all $v\in\cS_2^{d-1}$, $\inr{\tilde X_1 - \mu^*, v}$ has the same distribution as $\inr{\tilde X_1 - \mu^*, e_1}$ where $e_1=(1,0,\ldots,0)\in\bR^d$. In that case, $\tilde X_1$ has the same distribution as $\mu^*+R U$  where $R$ is a real-valued random variable on $\bR_+$ independent of $U$ a random vector uniformly distributed  over $\cS_2^{d-1}$. In that case and for $K=N$, for all $v\in\cS_2^{d-1}$ and all $r\in\bR$, 
\begin{equation*}
 H_{N,K=N,v}(r)=H(r): = \bP[R \inr{U,e_1}\geq r] = \int_r^{+\infty}f(x)dx \mbox{ where } f:x\in\bR\to C_d \int_{|x|}^{+\infty}\left(1-\frac{x^2}{u^2}\right)^{\frac{d-3}{2}} d \bP_R(u),
 \end{equation*} $\bP_R$ is the probability distribution of $R$ and $C_d$ is a normalization constant which can be proved to satisfy $\sqrt{d}\leq C_d\leq 6 \sqrt{d}$ (see for instance, Chapter~4 in \cite{MR1335228}). In particular, it follows from the mean value theorem that for all $r\geq 0$, $H(r)\leq H(0)-\min_{0\leq x\leq r}f(x)r = 1/2 - f(r)r$. Therefore, Assumption~\ref{ass:function_H} holds in that case when there exists constants $c_0,c_1>0$ such that $f(c_0)\geq c_1$. Furthermore, we have
 \begin{equation*}
 f(c_0)\geq C_d \int_{c_0\sqrt{d}}^{+\infty}\left(1-\frac{c_0^2}{u^2}\right)^{\frac{d-3}{2}} d \bP_R(u)\geq \frac{\sqrt{d}}{2}\bP[R\geq c_0\sqrt{d}]
 \end{equation*} because $C_d\geq \sqrt{d}$ and for all $u\geq c_0\sqrt{d}$, $(1-(c_0/u)^2)^{(d-3)/2}\geq 1/2$. As a consequence, Assumption~\ref{ass:function_H} holds if there are some constants $c_0>0$ and $c_2>0$ such that $\bP[R\geq c_0\sqrt{d}]\geq c_2/\sqrt{d}$. This is for instance the case, when $R$ is distributed like $\norm{G}_2$ for $G\sim \cN(0,I_d)$ (in that case $\tilde X_1\sim \cN(\mu^*, I_d)$) because $\bP[\norm{G}_2\geq \bE\norm{G}_2/2]\geq 1/2$ and $\bE \norm{G}_2\geq \sqrt{d}/2$ by Borell-TIS inequality but as well when $R$ is the positive part of a Cauchy variable because $\int_{\sqrt{d}}^{+\infty}(1/(1+x^2)) dx\geq 1/(2\sqrt{d})$. As a consequence, Assumption~\ref{ass:function_H} has nothing to do with the existence of any moment and it may hold even when there is not a first moment and even for $K=N$.

 Another example where Assumption~\ref{ass:function_H} holds, that we will use in the following to obtain statistical bounds for the coordinate-wise median of means for the location problem is when $S=\{\pm e_j:j\in[d]\}$ and $\tilde X_1 = \mu^* + Z$ where $Z=(z_j)_{j=1}^d$ is random vector in $\bR^d$ with coordinates $z_1, \ldots, z_d$ having a symmetric around $0$ Cauchy distribution. In that case, $\tilde X_1$ does not have a first moment and $\mu^*$ is a location parameter as the center of symmetry of the distribution of $\tilde X_1$. We have for all $j\in [d]$, 
 \begin{equation*}
 H_{N,K=N,\pm e_j}(r) = \bP\left[\inr{\tilde X_1-\mu^*,\pm e_j}\geq r\right] = \bP[z_j\geq r]= \int_r^{+\infty}\frac{dx}{\pi(1+x^2)}\leq \frac{1}{2} - \frac{r}{\pi(1+r^2)}\leq \frac{1}{2} -  \frac{r}{2\pi}
 \end{equation*}for all $0<r\leq 1$. Therefore, Assumption~\ref{ass:function_H} holds in that case as well.


\subsection{Statistical bounds for $\hat \mu_S^f$ and $\hat \mu_S^g$} 
\label{sub:statistical_bounds_for_}
In this section, we obtain estimation bounds w.r.t. $\norm{\cdot}_S$ for $\hat \mu_S^f$ and $\hat \mu_S^g$ in the adversarial contamination model with either the $L_2$ moment Assumption~\ref{assum:first} or the regularity at $0$ Assumption~\ref{ass:function_H}.
\paragraph{Estimation properties of $\hat \mu_S^f$ and $\hat \mu_S^g$ under Assumption~\ref{assum:first}.} 
\label{sub:estimation_properties_of}
In this section, we  obtain high probability estimation upper bounds satisfied by $\hat \mu_S^f$ and $\hat \mu_S^g$ w.r.t. $\norm{\cdot}_S$ in the adversarial contamination and heavy-tailed inlier model. The rate of convergence is given by the quantity
\begin{equation}\label{eq:main_rate}
r_S^*=\max\left( \frac{64}{\sqrt{N}} \bE \norm{\frac{1}{\sqrt{N}}\sum_{i\in [N]}\eps_i (\tilde X_i-\mu^*)}_{S}, \sup_{v\in  S} \norm{\Sigma^{1/2} v}_2\sqrt{\frac{64 K}{N}}\right).
\end{equation}

The key metric property satisfied by the two Fenchel-Legendre transforms $f^*_S$ and $g^*_S$ in the adversarial contamination and heavy-tailed inlier model is the following isomorphic result. 
\begin{Lemma}\label{lemm:VecteurLemme}
Grant Assumption~\ref{assum:first} and Assumption~\ref{ass:moment}. Let $S$ be a symmetric subset of $\bR^d$. Assume that $|\cO|<K/16$. With probability at least $1- \exp(-K/512)$, for all $\mu\in \bR^d$, $\left|g^*_{S}(\mu) - \norm{\mu - \mu^*}_S\right|\leq g^*_{S}(\mu^*)\leq r_S^*$ and $\left|f^*_{S}(\mu) - \norm{\mu - \mu^*}_S\right|\leq f^*_{S}(\mu) \leq r_S^*$.
\end{Lemma}

Lemma~\ref{lemm:VecteurLemme} shows that if $\norm{\mu - \mu^*}_S\geq 2 r^*_S$ then $\norm{\mu - \mu^*}_S\leq g^*_{S}(\mu)\leq 2\norm{\mu - \mu^*}_S$ and the same holds for $f^*_{S}$. It means that both $g^*_{S}$ and  $f^*_{S}$ are two convex functions equivalent (up to absolute constants) to $\mu\to\norm{\mu - \mu^*}_S$ on $\bR^d\backslash (2r^*_S)B_S$, where $B_S$ is the unit ball associated with $\norm{\cdot}_S$ and,  on $(2r^*_S)B_S$, they are both smaller than $2r^*_S$. Hence, both $g^*_{S}(\cdot-\mu^*)$ and  $f^*_{S}(\cdot-\mu^*)$ provide a good approximation of the metric space $(\bR^d, \norm{\cdot}_S)$. In particular, any minimum of $g^*_{S}$ and  $f^*_{S}$ will be close (up to $r^*_S$) to a minimum of  $\mu\to\norm{\mu - \mu^*}_S$ which is $\mu^*$. This explains the statistical properties of $\hat \mu_S^f$ and $\hat \mu_S^g$: from Lemma~\ref{lemm:VecteurLemme},
\begin{equation*}
  \norm{\hat \mu^f_S - \mu^*}_S \leq f^*_S(\hat \mu^f_S) + f^*_S(\mu^*)\leq 2f^*_S(\mu^*)\leq 2 r^*_S
  \end{equation*}  and the same holds for $\hat \mu^g_S$. This leads to the following result.

\begin{Theorem}\label{theo:main_convex_progrm}
Grant Assumption~\ref{assum:first} and Assumption~\ref{ass:moment}. Let $S$ be a symmetric subset of $\bR^d$ and $r_S^*$ be defined in \eqref{eq:main_rate}. For all $K> 16 |\cO|$, with probability at least $1- \exp(-K/512)$,
\begin{equation*}
\norm{\hat \mu_S^f-\mu^*}_S\leq 2r_S^* \mbox{ and } \norm{\hat \mu_S^g-\mu^*}_S\leq 2r_S^*.
\end{equation*}
\end{Theorem}

The rate $r^*_S$ obtained in Theorem~\ref{theo:main_convex_progrm} can be split into two terms: the complexity term given by the Rademacher complexity and a deviation term exhibiting the weak variance term as in the Gaussian case. Compare with Theorem~\ref{theo:lugosi_mendelson_1} from \cite{MR4026610}, this result shows that the Gaussian mean width term appearing in Theorem~\ref{theo:lugosi_mendelson_1} is actually not necessary, it also shows that this improved rate can be obtained by a procedure solution to a convex program and that it can also handle  adversarial corruption. When $S=B_2^d$, we recover the classical subgaussian rate because in that case the Rademacher complexity term in $r_S^*$ is less or equal to $\sqrt{\Tr(\Sigma)}$ \cite{MR2329442}. In particular, since $\hat \mu^g_S$ is the minmax MOM estimator in that case, we recover the main result from \cite{lerasle2019monk}. 

\paragraph{Estimation properties of $\hat \mu_S^g$ under Assumption~\ref{ass:function_H}.} 
\label{sub:estimation_properties_of}
In this section, we consider some cases where a first moment may not exist; in that case, $\mu^*$ is a location parameter so that  Assumption~\ref{ass:function_H} holds. The rate of convergence we obtain in that case is given by
\begin{equation}\label{eq:r_diamond}
r^\diamond = \frac{C_0}{c_1}\left(\sqrt{\frac{d+1}{N}} + \sqrt{\frac{u}{N}}\right) + \frac{|\cO|}{c_1 \sqrt{K N}}
\end{equation}where $c_1$ is the absolute constant from Assumption~\ref{ass:function_H}, $C_0$ the absolute constant from \eqref{eq:VC_concentration} and $u>0$ a confidence parameter. 

The following result is an isomorphic result satisfied by the Fenchel-Legendre transforms  $g^*_S$ under Assumption~\ref{ass:function_H}. It is similar to the one of Lemma~\ref{lemm:VecteurLemme} but with the rate $r^\diamond$. 
\begin{Lemma}\label{lemm:VecteurLemme_under_ass:function_H} Let $S$ be a symmetric subset of $\bR^d$.
Grant Assumption~\ref{assum:first} and Assumption~\ref{ass:function_H} for some $K\in[N]$. Let $u>0$. Assume that $C_0\left(\sqrt{(d+1)/K} + \sqrt{u/K}\right) + |\cO|/K\leq c_0 c_1$. With probability at least $1- \exp(-u)$, for all $\mu\in \bR^d$, $\left|g^*_{S}(\mu) - \norm{\mu - \mu^*}_S\right|\leq r^\diamond$.
\end{Lemma}
As explained below Lemma~\ref{lemm:VecteurLemme}, a result such as Lemma~\ref{lemm:VecteurLemme_under_ass:function_H} may be used to upper bound the $\norm{\cdot}_S$ distance between $\hat \mu^g_S$, a minimum of $g^*_{S}$, and $\mu^*$, a minimum of $\mu\to\norm{\mu - \mu^*}_S$. This yields to the following result.

\begin{Theorem}\label{theo:main_2_regular_zero} Let $S$ be a symmetric subset of $\bR^d$.
Grant Assumption~\ref{assum:first} and Assumption~\ref{ass:function_H} for some $K\in[N]$.  Let $u>0$ and assume that $C_0\left(\sqrt{(d+1)/K} + \sqrt{u/K}\right) + |\cO|/K\leq c_0 c_1$.  With probability at least $1- \exp(-u)$, $\norm{\hat \mu_S^g-\mu^*}_S\leq 2r^\diamond$ where $r^\diamond$ is defined in \eqref{eq:r_diamond}.
\end{Theorem}

Unlike Theorem~\ref{theo:main_convex_progrm}, Theorem~\ref{theo:main_2_regular_zero} may hold even when there is not a first moment. The result from Theorem~\ref{theo:main_2_regular_zero} hold for all $0<u\lesssim K$ whereas Theorem~\ref{theo:main_convex_progrm} holds only for $u=K$ (even though one may use a Lepski's adaptive scheme to chose adaptively $K$). The price for adversarial corruption in \eqref{eq:r_diamond} is between $|\cO|/N$ (for $K \sim N$) and $\sqrt{|\cO|/N}$ (for $K\sim |\cO|$). It therefore depends on the choice of $K$ for which Assumption~\ref{ass:function_H} holds. As shown after Assumption~\ref{ass:function_H} for spherically symmetric random variables one can take $K=N$ and so the best possible price $|\cO|/N$ for adversarial corruption may be achieved even when a first moment does not exist. If one needs some averaging effect so that  Theorem~\ref{theo:main_2_regular_zero} holds, then one should take $K$ as small as possible that is $K\sim|\cO|$ and then  $\sqrt{|\cO|/N}$ will be the price for adversarial corruption as in the $L_2$ case described in Theorem~\ref{theo:main_2_regular_zero}.

\paragraph{Subgaussian rates under weak or no moment assumption.} 
\label{par:subgaussian_rates_under_weak_moment_assumptions_}
It is possible to recover (up to absolute constants) the subgaussian rate \eqref{eq:subgaussian_rate_S} in Theorem~\ref{theo:main_convex_progrm} for $K\sim \log(1/\delta)$ when the Rademacher complexity term from \eqref{eq:main_rate} and the Gaussian mean width from \eqref{eq:subgaussian_rate_S} satisfy
\begin{equation}\label{eq:condition_subgauss_rate}
\bE \norm{\frac{1}{\sqrt{N}}\sum_{i\in [N]}\eps_i (\tilde X_i-\mu^*)}_{S} \lesssim \ell^*\left(\Sigma^{1/2}S\right).
\end{equation}Such a result (i.e. Rademacher complexity is smaller than the Gaussian mean width up to constant) depends on the set $S$ and the number of moments granted on the $\tilde X_i$'s as well as the sample size. It obviously holds when the $\tilde X_i$'s are i.i.d. $\cN(\mu^*, \Sigma)$, so that we recover the deviation-minimax optimal subgaussian rate \eqref{eq:subgaussian_rate_S} in that case. It is also true when the $\tilde X_i$'s are subgaussian vectors. There are other situations under weaker moment assumption where \eqref{eq:condition_subgauss_rate} holds.

For instance, when $S=B_2^d$, \eqref{eq:condition_subgauss_rate} holds under only a $L_2$-moment assumption (see \cite{MR2329442}). It also holds for $S=B_1^d$ when the $\tilde X_i$'s are isotropic with  coordinates having $\log d$ subgaussian moments (i.e. $\norm{\inr{\tilde X_i,e_j}}_{L_p}\leq L \sqrt{p}$ for all $1\leq p\leq \log d$ and coordinate $j\in[d]$) and $N\gtrsim \log d$. Together with \eqref{eq:coordinate_wise_MOM} and Theorem~\ref{theo:main_convex_progrm}, this implies that the coordinate-wise MOM is a subgaussian estimator of the mean under a $\log d$ subgaussian moment assumption.  Upper bounds such as \eqref{eq:condition_subgauss_rate} have been extended in  \cite{MR3645129} to general unconditional norms. 

It is also possible to recover the subgaussian rate \eqref{eq:subgaussian_rate_S}  in situations where there is not even a first moment thanks to Theorem~\ref{theo:main_2_regular_zero}. Indeed, for the case $S=B_1^d$ and $\tilde X_1=\mu^* + Z$ where $Z=(z_j)_{j=1}^d$ has symmetric around $0$ Cauchy distributed coordinates, we showed that Assumption~\ref{ass:function_H} holds for $K=N$ and that $\hat \mu_S^g$ is the coordinate-wise median (here $K=N$) in \eqref{eq:coordinate_wise_MOM}. It follows from Theorem~\ref{theo:main_2_regular_zero} that,  when $d\lesssim N$ and $|\cO|\lesssim N$ then for all $d\leq u\lesssim N$, with probability at least $1-\exp(-u)$,
\begin{equation}\label{eq:estim_prop_coordinate_wise_med}
\norm{\hat\mu^g_S - \mu^*}_\infty\leq 2C_0\left(\sqrt{\frac{d+1}{N}} + \sqrt{\frac{u}{N}}\right) + \frac{2\pi |\cO|}{N}
\end{equation}which is the deviation-minimiax optimal subgaussian rate \eqref{eq:subgaussian_rate_S} we would have gotten if the $\tilde X_i$ were i.i.d. isotropic Gaussian vectors centered in $\mu^*$ corrupted by $|\cO|$ adversarial outliers (up to absolute constants). But here, \eqref{eq:estim_prop_coordinate_wise_med} is obtained without the existence of a first moment. Moreover, in \eqref{eq:estim_prop_coordinate_wise_med}, the number of outliers is allowed to be proportional to $N$ and the price for adversarial corruption is of the order of $|\cO|/N$ which is the same price we have to pay when inlier have a Gaussian distribution -- this differs from the $\sqrt{|\cO|/N}$ information theoretical lower bound that has been obtained for some non-symmetric inlier.  Furthermore, the computational cost of the coordinate-wise MOM is $\cO(Nd)$ since the cost for computing the bucketed means is $\cO(Nd)$, the one of finding the median of $K$ numbers is $\cO(K)$ \cite{blum1973time}, it is therefore the same computational cost as the one of the empirical mean. It is therefore possible to achieve the same computational and statistical properties as the empirical mean in a setup where a first moment does not even exist.

\section{Proofs} 
\label{sec:proofs}
\textbf{Proof of Theorem~\ref{theo:minimax_gauss_S}.}
The minimax lower bound rate $r^*$ exhibits two quantities: one which is a \textit{complexity term} depending on the Gaussian mean width of $\Sigma^{1/2} S$ and a \textit{deviation term} depending on $\delta$. The two terms come from two arguments. We start with the deviation term.

Let $v_1\in\bR^d$ be such that $\norm{v_1}_S=1$. We consider two Gaussian measures on $\bR^{dN}$: $\bP_0=\cN(0,\Sigma)^{\otimes N}$ and $\bP_1 = \cN(3 r^*v_1, \Sigma)^{\otimes N}$. They are the distributions of a sample of $N$ i.i.d. Gaussian vectors in $\bR^d$ with the same covariance matrix $\Sigma$ and the first one with mean $0$ and the second one with mean $3r^* v_1$. We set 
$A_0 = (\hat \mu)^{-1}(B_S(0, r^*)) = \{(x_1, \ldots, x_N)\in\bR^{Nd}: \norm{\hat \mu(x_1, \ldots, x_N) }_S\leq r^*\}$ and $A_1=(\hat \mu)^{-1}(B_S(3r^*v_1, r^*))$. It follows from the statistical properties of $\hat \mu$ that $\bP_0[A_0]\geq 1-\delta$ and $\bP_1[A_1]\geq 1 - \delta$. 

The key ingredient for the deviation lower bound term is a slightly generalization of Lemma~3.3 in \cite{lecue2013learning} which is  based on a version of the Gaussian shift Theorem from \cite{MR1652329}.

\begin{Lemma} \label{lemma:distance-vs-measure}
Let $t\mapsto\Phi(t)=\bP(g \leq t)$ be the cumulative distribution function of a standard gaussian random variable on $\R$.
Let $\Sigma_0\succeq 0$ be in $\bR^{(Nd)\times (Nd)}$ and $u,v \in \R^{dN}$. Let two gaussian measures $\nu_u\sim \cN(u,\Sigma_0)$ and $\nu_v\sim\cN(v,\Sigma_0)$ on $\bR^{Nd}$. If $A \subset \R^{dN}$ is measurable, then
\begin{equation}\label{eq:shift_theorem}
\nu_v(A) \geq 1-\Phi\big(\Phi^{-1}(1-\nu_u(A))+\|\Sigma_0^{-1/2}(u-v)\|_{2}\big)
\end{equation}where $\Sigma_0^{-1/2}$ is the square root of the pseudo-inverse of $\Sigma_0$.
\end{Lemma}
\textbf{Proof of Lemma~\ref{lemma:distance-vs-measure}.} When $\Sigma_0 = I_{Nd}$, Lemma~\ref{lemma:distance-vs-measure} is exactly Lemma~3.3 in \cite{lecue2013learning} for $\sigma=1$. To prove Lemma~\ref{lemma:distance-vs-measure}, we observe that $\nu_v(A) = \bP[G+\Sigma_0^{-1/2}v\in B]$ where $B=\Sigma_0^{-1/2}A$ and $G$ is a standard Gaussian variable in ${\rm Im}(\Sigma_0)$. Hence, it follows from Lemma~3.3 in \cite{lecue2013learning} that 
\begin{equation*}
\bP[G+\Sigma_0^{-1/2}v\in B] \geq 1-\Phi\big(\Phi^{-1}(1-\bP[G+\Sigma_0^{-1/2}u\in B])+\|\Sigma_0^{-1/2}(u-v)\|_{\ell_2^N}\big)
\end{equation*}which is exactly \eqref{eq:shift_theorem}. {\mbox{}\nolinebreak\hfill\rule{2mm}{2mm}\par\medbreak}

It follows from Lemma~\ref{lemma:distance-vs-measure} that 
\begin{equation}\label{eq:shift_ineq}
\bP_1[A_0]\geq 1 - \Phi\left[\Phi^{-1}(1-\bP_0[A_0])+\norm{\Sigma_0^{-1/2}(0 - (3r^*v_1, \ldots, 3r^*v_1))}_2\right].
\end{equation}Moreover, we have $\Phi^{-1}(1-\bP_0[A_0])\leq \Phi^{-1}(\delta)$ (because $1-\bP_0[A_0]\leq\delta$) and 
\begin{equation}\label{eq:big_var}
\norm{\Sigma_0^{-1/2}(0 - (3r^*v_1, \ldots, 3r^*v_1))}_2 = 3r^*\sqrt{N}\norm{\Sigma^{-1/2}v_1}_2.
\end{equation}As a consequence, if $3r^*\sqrt{N} \norm{\Sigma^{-1/2}v_1}_2\leq -\Phi^{-1}(\delta)$ then, in \eqref{eq:shift_ineq}, we get $\bP_1[A_0]\geq 1 - \Phi[0]\geq 1/2$ which is not possible because $\bP_1[A_1]\geq 1-\delta>3/4$ and $A_1\cap A_0=\emptyset$. As a consequence, we necessarily have $3r^*\sqrt{N}\geq (-\Phi^{-1}(\delta)) \norm{\Sigma^{-1/2}v_1}_2^{-1}$. The later holds for any $v_1\in\bR^d$ such that $\norm{v_1}_S=1$ hence $3r^*\sqrt{N}\geq (-\Phi^{-1}(\delta))[1/\inf_{\norm{v}_S=1}\norm{\Sigma^{-1/2}v}_2]$. It also follows from the bound on the Mill's ratio from \cite{MR0079844} (here we use that for all $x \geq0$, $\Phi(-x)\geq 2 \varphi(x)/{\sqrt{4+x^2}+x}$ where $\varphi$ is the standard Gaussian density function) that for all $0<\delta<1/4$, $-\Phi^{-1}(\delta)\geq 1/4\sqrt{\log(1/\delta)}$. This shows that
\begin{equation}\label{eq:final_1_before_dual}
r^*\geq \frac{1}{12}\sqrt{\frac{\log(1/\delta)}{N}}\frac{1}{\inf_{\norm{v}_S=1}\norm{\Sigma^{-1/2}v}_2}.
\end{equation}To conclude on the deviation term, we use the following duality argument.

\begin{Lemma}\label{lem:dual}
 Let $A\in\bR^{d\times d}$ be a symmetric and invertible matrix. Let $\norm{\cdot}$ be a norm and its dual norm $\norm{\cdot}^*$ on $\bR^d$. Let $S$ be a symmetric subset of $\bR^d$ such that ${\rm span}(S) = \bR^d$. We have
 \begin{equation*}
 \frac{1}{\inf_{\norm{v}_S=1}\norm{A^{-1}v}} \geq \sup_{w\in S} \norm{A w}^*.
 \end{equation*}
 \end{Lemma} 
 \textbf{Proof of Lemma~\ref{lem:dual}.} Let $v$ be such that $\norm{v}_S=1$ and  $w\in S$. We have $|\inr{v,w}|\leq 1$ and so $|\inr{A^{-1}v/\norm{A^{-1}v}, Aw}|\leq 1 / \norm{A^{-1}v}$. The later holds for all $v$ such that $\norm{v}_S=1$ and  $\{A^{-1}v/\norm{A^{-1}v}:\norm{v}_S=1\}$ is the unit sphere of $\norm{\cdot}$. Hence, we conclude by taking the sup over $v$ such that $\norm{v}_S=1$ and $w\in S$. {\mbox{}\nolinebreak\hfill\rule{2mm}{2mm}\par\medbreak}

 It follows from \eqref{eq:final_1_before_dual} and Lemma~\ref{lem:dual} for $\norm{\cdot} = \norm{\cdot}_2$ and $A=\Sigma^{1/2}$ that 
 \begin{equation}\label{eq:deviation_bound}
r^*\geq \frac{1}{12}\sqrt{\frac{\log(1/\delta)}{N}}\sup_{w\in S}\norm{\Sigma^{1/2}w}_2.
 \end{equation}

Let us now turn to the second part of the lower bound; the one coming from the complexity of the problem (here, it is the Gaussian mean width of $\Sigma^{1/2}S$).  We know that $\hat\mu$ is an estimator such that for all $\mu\in\bR^d$, $\bP_\mu^N\left[\norm{\hat \mu - \mu}_S\leq r^*\right]\geq 1-\delta$ which is equivalent to say that 
\begin{equation}\label{eq:risk_minmax}
 \delta \geq \sup_{\mu\in\bR^d} \bE_\mu^N \phi\left(\frac{\norm{\hat\mu - \mu}_S}{r^*}\right)
 \end{equation} where we set $\phi:t\in\bR \to I(t>1)$ and $\bE^N_\mu$ is the expectation with respect to $X_1, \ldots X_N \ove{\sim}{i.i.d.}\cN(\mu, \Sigma)$. Next, we consider a Gaussian distribution $\gamma$ over the set of parameters $\mu\in\bR^d$: for $s>0$, we assume that $\mu\sim \cN(0, s \Sigma)$. It follows from \eqref{eq:risk_minmax} that 
 \begin{equation}\label{eq:apply_anderson}
  \delta\geq \int_{\mu\in\bR^d} \bE_\mu^N \phi\left(\frac{\norm{\hat\mu - \mu}_S}{r^*}\right) \gamma(\mu)d\mu = \bE\left[\bE\left[ \phi\left(\frac{\norm{\hat\mu(X_1, \ldots, X_N) - \mu}_S}{r^*}\right)  | X_1, \ldots, X_N\right] \right]. 
  \end{equation} In other words, we lower bound the minmax risk by a Bayesian risk. We now use Anderson's lemma to lower bound the Bayesian risk appearing in \eqref{eq:apply_anderson}. We first recall Anderson's Lemma.

  \begin{Theorem}[Anderson's Lemma] Let $\Gamma$ be a semi-definite $d\times d$ matrix and $Z\sim\cN(0,\Gamma)$. Let $w:\bR^d\to \bR$ be such that all its level sets (i.e. $\{x\in\bR^d:w(x)\leq c\}$ for $c\in\bR$) are convex and symmetric around the origin. Then for all $x\in\bR^d$, $\bE w(Z+x)\geq \bE w(Z)$.
  \end{Theorem}

We remark that $\mu - \bE[\mu|X_1, \ldots, X_N]$ is distributed according to $\cN(0, (s/(1+Ns) \Sigma))$ conditionally to $X_1, \ldots,X_N$. Therefore, applying Anderson's Lemma conditionally to $X_1, \ldots, X_N$, we obtain in \eqref{eq:apply_anderson} that
\begin{equation*}
 \delta\geq \bE \left[ \phi\left(\frac{\norm{\bE[\mu|X_1, \ldots, X_N] - \mu}_S}{r^*}\right)\right] = \bP\left[\norm{\Sigma^{1/2} G}_S\geq \sqrt{\frac{1+Ns}{s}} r^*\right]
 \end{equation*}where $G\sim \cN(0, I_d)$. This result is true for all $s>0$ so taking $s\uparrow+\infty$, we obtain
 \begin{equation*}
 \delta\geq \bP\left[\norm{\Sigma^{1/2}G}_S\geq \sqrt{N} r^*\right].
 \end{equation*}Using Borell-TIS's inequality (Theorem~7.1 in \cite{Led01} or pages 56-57 in \cite{MR3184689}), we know that  with probability at least $4/5$, $\norm{\Sigma^{1/2}G}_S\geq \bE \norm{\Sigma^{1/2}G}_S-  \sigma_S\sqrt{2 \log(5/4)}$ where we set $\sigma_S = \sup_{\norm{v}_S=1}\norm{\Sigma^{1/2}v}_2$. As a consequence, for $\delta=1/4$, we necessarily have $\sqrt{N} r^*\geq \bE \norm{\Sigma^{1/2}G}_S -  \sigma_S\sqrt{2\log(5/4)}$ and so $\sqrt{N} r^*\geq (1/2)\bE\norm{\Sigma^{1/2}G}_S$ when $\bE\norm{\Sigma^{1/2}G}_S\geq 2\sigma_S\sqrt{2\log(5/4)}$. Finally, when $\bE\norm{\Sigma^{1/2}G}_S< 2\sigma_S\sqrt{2\log(5/4)}$, we know from \eqref{eq:deviation_bound} for $\delta=1/4$ that
 \begin{equation*}
 r^*\geq \frac{1}{12}\sqrt{\frac{\log4}{N}}\sigma_S\geq \frac{1}{24}\sqrt{\frac{\log 2}{\log(5/4)}}\frac{\bE\norm{\Sigma^{1/2}G}_S}{\sqrt{N}}.
 \end{equation*}{\mbox{}\nolinebreak\hfill\rule{2mm}{2mm}\par\medbreak}

\paragraph{Proof of Theorem~\ref{theo:minimax_gauss}.}
\label{par:proof_of_theorem_minimax_gauss} 
Theorem~\ref{theo:minimax_gauss} follows from Theorem~\ref{theo:minimax_gauss_S} and the following lower bound on $\bE \norm{\Sigma^{1/2}G}_{B_2^d}$. We have from Borell-TIS's inequality that 
  \begin{align*}
  &\bE \norm{\Sigma^{1/2}G}_2^2 - \left(\bE \norm{\Sigma^{1/2}G}_2\right)^2 = \bE\left(\norm{\Sigma^{1/2}G}_2-\bE\norm{\Sigma^{1/2}G}_2 \right)^2\\
  & = \int_0^\infty \bP\left[ \left|\norm{\Sigma^{1/2}G}_2-\bE\norm{\Sigma^{1/2}G}_2 \right|\geq\sqrt{t} \right]dt\leq 2 \sigma^2_{B_2^d}
  \end{align*}where $\sigma^2_{B_2^d} = \sup_{\norm{v}_2=1}\norm{\Sigma^{1/2}v}_2^2 = \norm{\Sigma}_{op}$. Since $\bE \norm{\Sigma^{1/2}G}_2^2 = {\rm Tr}(\Sigma)$, we have $\left(\bE \norm{\Sigma^{1/2}G}_2\right)^2\geq {\rm Tr}(\Sigma) - 2\norm{\Sigma}_{op}$. Therefore, $\bE \norm{\Sigma^{1/2}G}_2\geq \sqrt{{\rm Tr}(\Sigma)/2}$ when ${\rm Tr}(\Sigma) \geq  4\norm{\Sigma}_{op}$ and when ${\rm Tr}(\Sigma) <  4\norm{\Sigma}_{op}$, we use the lower bound from \eqref{eq:deviation_bound} and an argument similar to the one appearing in the end of the proof of Theorem~\ref{theo:minimax_gauss_S} to get the result.
  {\mbox{}\nolinebreak\hfill\rule{2mm}{2mm}\par\medbreak}

\paragraph{Proof of Lemma~\ref{lemm:VecteurLemme}.} 
\label{par:proof_of_lemma_lemm:vecteurlemme}
We first prove the result for the $g^*_S$ function. The one for the $f^*_S$ is similar up to constants and will be sketched after. The proof of Lemma~\ref{lemm:VecteurLemme} for the $g^*_S$ function is a corollary of the general fact which holds under only Assumption~\ref{assum:first}. Let $u>0$ be a confidence parameter and define $R^*_S$ such that
\begin{equation}\label{eq:demo_lem1_fixed_point}
\frac{4}{\sqrt{N}R^*_S}\bE \norm{\frac{1}{\sqrt{N}}\sum_{i\in [N]}\eps_i (\tilde X_i-\mu)}_{S}  +  \sqrt{\frac{2u}{K}} + \sup_{v\in S}H_{N,K,v}\left(\frac{R^*_S}{2}\sqrt{\frac{N}{K}}\right) + \frac{|\cO|}{K} < \frac{1}{2}.
\end{equation}Let us show that with large probability for all $\mu\in\bR^d$, $\left|g^*_S(\mu) - \norm{\mu - \mu^*}_S\right|\leq R^*_S$.

We have for all $\mu\in\bR^d$, 
\begin{equation}\label{eq:demo_lem1_first}
\left|g^*_S(\mu) - \norm{\mu - \mu^*}_S\right| = \left|\sup_{v\in S}\left(\inr{\mu,v} - g(v)\right) - \sup_{v\in S}\inr{v, \mu-\mu^*}\right|\leq \sup_{v\in S}\left|\inr{\mu^*,v} - g(v)\right|=g^*_S(\mu^*)
\end{equation}where we used that $S$ is symmetric and $g$ is odd. It only remains to show that $g^*_S(\mu^*)\leq R^*_S$ with large probability. To that end, it is enough to prove that, with large probability, for all $v\in S$, 
 \begin{equation}\label{eq:result_bucket_mean}
 \sum_{k\in[K]} I(\inr{\bar X_k - \mu^*, v}> R^*_S)< \frac{K}{2}.
 \end{equation}

We use the notation introduced in Assumption~\ref{assum:first} and we consider ${\overline{\tilde X}}_k= |B_k|^{-1}\sum_{i\in B_k} \tilde X_i$ for $k\in[K]$ which are the  $K$ bucketed means constructed on the $N$ independent vectors $\tilde X_i, i\in[N]$ before contamination (whereas $\bar X_k$ are the ones constructed after contamination). We also set  $\cK=\{k\in[K]:B_k\cap\cO=\emptyset\}$ the indices of the non corrupted blocks. We have
 \begin{align}\label{eq:result_bucket_mean_2}
 \notag &\sum_{k\in[K]} I(\inr{\bar X_k - \mu^*, v}> R^*_S) = \sum_{k\in\cK} I(\inr{\bar X_k - \mu^*, v}> R^*_S) + \sum_{k\notin\cK} I(\inr{\bar X_k - \mu^*, v}> R^*_S)\\
 &\leq \sum_{k\in[K]} I(\inr{{\overline{\tilde X}}_k - \mu^*, v}> R^*_S)  + |\cO|.
 \end{align}

It only remains to show that with probability at least $1-\exp(-u)$, for all $v\in  S$, 
 \begin{equation*}
 \sum_{k\in[K]} I(\inr{{\overline{\tilde X}}_k - \mu^*, v}> R^*_S)\leq \frac{4K}{\sqrt{N}R^*_S}\bE \norm{\frac{1}{\sqrt{N}}\sum_{i\in [N]}\eps_i (\tilde X_i-\mu^*)}_{S}  +  \sqrt{2uK} + K\sup_{v\in S}H_{N,K,v}\left(\frac{R^*_S}{2}\sqrt{\frac{N}{K}}\right).
 \end{equation*}

 We define $\phi(t) = 0 $ if $t\leq1/2$, $\phi(t) = 2(t-1/2)$ if $1/2\leq t\leq 1$ and $\phi(t) = 1$ if $t\geq1$. We have $I(t\geq1)\leq \phi(t)\leq I(t\geq1/2)$ for all $t\in\bR$ and so
  \begin{align*}
 &\sum_{k\in [K]} I(\inr{{\overline{\tilde X}}_k - \mu^*, v}> R^*_S)\\
 &\leq \sum_{k\in[K]} I(\inr{{\overline{\tilde X}}_k - \mu^*, v}> R^*_S) - \bP[\inr{{\overline{\tilde X}}_k - \mu^*, v}> R^*_S/2] + \bP[\inr{{\overline{\tilde X}}_k - \mu^*, v}> R^*_S/2]\\
 &\leq \sum_{k\in[K]}\phi\left(\frac{\inr{{\overline{\tilde X}}_k - \mu^*, v}}{R^*_S}\right) - \bE \phi\left(\frac{\inr{{\overline{\tilde X}}_k - \mu^*, v}}{R^*_S}\right) + \bP[\inr{{\overline{\tilde X}}_k - \mu^*, v}> R^*_S/2]\\
& \leq \sup_{v\in  S}\left(\sum_{k\in[K]}\phi\left(\frac{\inr{{\overline{\tilde X}}_k - \mu^*, v}}{R^*_S}\right) - \bE \phi\left(\frac{\inr{{\overline{\tilde X}}_k - \mu^*, v}}{R^*_S}\right) \right)+  K \sup_{v\in S}H_{N,K,v}\left(\frac{R^*_S}{2}\sqrt{\frac{N}{K}}\right).
 \end{align*}

 Next, we use several tools from empirical process theory and in particular, for a symmetrization argument, we consider a family of $N$ independent Rademacher variables $(\eps_i)_{i=1}^N$ independent of the $(\tilde X_i)_{i=1}^N$.  In \textit{(bdi)} below, we use  the bounded difference inequality (Theorem~6.2 in \cite{MR3185193}). In \textit{(sa-cp)}, we use the symmetrization argument and the contraction principle (Chapter~4 in \cite{MR2814399}) -- we refer to the supplementary material of \cite{LMSL} for more details. We have,  with probability at least $1-\exp(-u)$, 
 \begin{align*}
  &\sup_{v\in  S}\left(\sum_{k\in[K]}\phi\left(\frac{\inr{{\overline{\tilde X}}_k - \mu^*, v}}{R^*_S}\right) - \bE \phi\left(\frac{\inr{{\overline{\tilde X}}_k - \mu^*, v}}{R^*_S}\right) \right)\\
  & \ove{\leq}{(bdi)} \bE \sup_{v\in  S}\left(\sum_{k\in[K]}\phi\left(\frac{\inr{{\overline{\tilde X}}_k - \mu^*, v}}{R^*_S}\right) - \bE \phi\left(\frac{\inr{{\overline{\tilde X}}_k - \mu^*, v}}{R^*_S}\right) \right)+ \sqrt{2uK}\\
  &\ove{\leq}{(sa-cp)} \frac{4K}{NR^*_S} \bE \sup_{v\in  S} \inr{v, \sum_{i\in [N]}\eps_i (\tilde X_i-\mu^*)} + \sqrt{2uK}\\
  & =  \frac{4K}{\sqrt{N}R^*_S}\bE \norm{\frac{1}{\sqrt{N}}\sum_{i\in [N]}\eps_i (\tilde X_i-\mu^*)}_{S}  + \sqrt{2uK}.
 \end{align*}

We therefore showed that under Assumption~\ref{assum:first}, with probability at least $1-\exp(-u)$, for all $\mu\in\bR^d$, $\left|g^*_S(\mu) - \norm{\mu - \mu^*}_S\right|\leq R^*_S$.

Now, if Assumption~\ref{ass:moment} holds then for all $v\in S$, we have from Markov's inequality that
 \begin{align*}
 &H_{N,K,v}\left(\frac{R^*_S}{2}\sqrt{\frac{N}{K}}\right)\leq \frac{\bE \inr{{\overline{\tilde X}}_k - \mu, v}^2}{(r^*_S/2)^2} = \frac{4Kv^\top \Sigma v}{N(r^*_S)^2}\leq  \frac{4K\sup_{v\in  S} \norm{\Sigma^{1/2} v}_2^2}{N(r^*_S)^2} \leq \frac{1}{8}
 \end{align*}and therefore \eqref{eq:demo_lem1_fixed_point} holds for $R^*_S = r^*_S$ when $|\cO|< K/8$ and $u=K/128$. This proves the result of Lemma~\ref{lemm:VecteurLemme} for $g^*_S$ under Assumption~\ref{ass:moment}.

 Finally, for the function $f^*_S$ one needs to control the average of the $K/2$ inter-quartiles. One way to do it is to control the value of all elements $\inr{\bar X_k - \mu^*, v}$ in the inter-quartiles interval. This can be done by defining an $R^*_S$ similar to the one in \eqref{eq:demo_lem1_fixed_point} but where the right-hand side value $1/2$  is replaced by $1/4$ in \eqref{eq:demo_lem1_fixed_point}. This only modifies the absolute constants which are the one used in  Lemma~\ref{lemm:VecteurLemme}.
 {\mbox{}\nolinebreak\hfill\rule{2mm}{2mm}\par\medbreak}

\paragraph{Proof of Lemma~\ref{lemm:VecteurLemme_under_ass:function_H}.} 
\label{par:proof_of_lemma}
Unlike in Lemma~\ref{lemm:VecteurLemme} where we used the Rademacher complexities as a complexity measure, in this proof, the complexity measure we are using is the Vapnik and Chervonenkis (VC) dimension \cite{MR3408730,MR1719582} of a class $\cF$ of Boolean functions, i.e. of functions from $\bR^d$ to $\{0,1\}$ in our case. We recall that the Vapnik and Chervonenkis dimension of $\cF$, denoted by $VC(\cF)$, is the maximal integer $n$ such that there exists $x_1, \ldots, x_n\in\bR^d$ for which the set $\{(f(x_1), \ldots, f(x_n)):f\in\cF)\}$ is of maximal cardinality, that is of size $2^n$. The VC dimension of the set of all indicators of half affine spaces in $\bR^d$ is  $d+1$ (see Example~2.6.1 in \cite{MR1385671}). We also know (see, for instance, Chapter~3 in \cite{Kol11}) the following concentration bound: let $Y_1, \ldots, Y_n$ be independent random vectors in $\bR^d$, there exists an absolute constant $C_0$ such that for all $u>0$, with probability at least $1-\exp(-u)$,
\begin{equation}\label{eq:VC_concentration}
\sup_{f\in\cF}\left(\frac{1}{n}\sum_{i=1}^n f(Y_i)-\bE f(Y_i)\right)\leq C_0\left(\sqrt{\frac{VC(\cF)}{n}} + \sqrt{\frac{u}{n}}\right).
\end{equation}

Lemma~\ref{lemm:VecteurLemme_under_ass:function_H} is a corollary of a general result which holds under the only Assumption~\ref{assum:first}. This general result says that for all $u>0$, with probability at least $1- \exp(-u)$, for all $\mu\in \bR^d$, $\left|g^*_{S}(\mu) - \norm{\mu - \mu^*}_S\right|\leq R^\diamond$ where $R^\diamond$ is any point such that 
\begin{equation}\label{eq:key_fixed_point}
C_0\left(\sqrt{\frac{d+1}{K}}+\sqrt{\frac{u}{K}}\right) + \sup_{\norm{v}_2=1}H_{N,K,v}\left(R^\diamond\sqrt{\frac{N}{K}}\right) + \frac{|\cO|}{K}< \frac{1}{2}
\end{equation}where $C_0$ is the constant from \eqref{eq:VC_concentration}. In particular, when Assumption~\ref{ass:function_H} holds then one can check that \eqref{eq:key_fixed_point} holds for $R^\diamond=r^\diamond$ when $r^\diamond\leq c_0 $ proving the result of Lemma~\ref{lemm:VecteurLemme_under_ass:function_H}. It only remains to show the general result. To that end we follow the same strategy as in the proof of Lemma~\ref{lemm:VecteurLemme} up to \eqref{eq:result_bucket_mean_2} (and with $R^*_S$ replaced by $R^\diamond$). From that point, we use \eqref{eq:VC_concentration} and the VC dimension of the set of affine half spaces to get that with probability at least $1-\exp(-u)$, for all $v\in S$, 
\begin{align*}
\sum_{k\in[K]} I(\inr{{\overline{\tilde X}}_k - \mu^*, v}> R^\diamond) \leq H_{N,K,v}\left(R^\diamond\sqrt{\frac{N}{K}}\right) + C_0\left(\sqrt{\frac{d+1}{N/K}} + \sqrt{\frac{u}{N/K}}\right)
\end{align*}and so by definition of $R^\diamond$, on the same event, for all  $v\in S$, $\sum_{k\in[K]} I(\inr{{{\bar X}}_k - \mu^*, v}> R^\diamond) <1/2$. This concludes the proof.
 {\mbox{}\nolinebreak\hfill\rule{2mm}{2mm}\par\medbreak}






\begin{footnotesize}
\bibliographystyle{plain}
\bibliography{biblio}
\end{footnotesize}


\end{document}

%% file: packages_notes.tex
\usepackage{graphicx}
\usepackage{color}
\usepackage{amsmath}
\usepackage{amssymb}
\usepackage{amscd}
\usepackage{bbm}

\usepackage{tikz}
\newlength\Colsep
\usetikzlibrary{patterns}

\usepackage{hyperref}
\hypersetup{
    bookmarks=true,         
    colorlinks=true,       
    linkcolor=red,          
    citecolor=green,        
    filecolor=magenta,      
    urlcolor=cyan           
}

\usepackage[utf8]{inputenc}
\usepackage{amsthm}
\usepackage{bbm} 
\usepackage[linesnumbered,lined,boxed,commentsnumbered]{algorithm2e}

\usepackage{marginnote}

\newtheorem{Theorem}{Theorem}
\newtheorem{Assumption}{Assumption}
\newtheorem{Definition}{Definition}
\newtheorem{Lemma}{Lemma}
\newtheorem{Proposition}{Proposition}


%% file: commandes_guillaume.tex
\newcommand{\inr}[1]{\bigl< #1 \bigr>}

\newcommand{\norm}[1]{\left\|#1\right\|}%

\newcommand\ove[2]{\overset{#2}{#1}\;}

\newcommand\eps{\epsilon}

\newcommand{\beginproof}{{\bf Proof. {\hspace{0.2cm}}}}
\def \endproof
{{\mbox{}\nolinebreak\hfill\rule{2mm}{2mm}\par\medbreak}}

\DeclareMathOperator*{\argmin}{argmin}

\def\ds1{\textrm{1\kern-0.25emI}} 

\newcommand \R{\mathbb{R}}

\newcommand \cF{{\cal F}}

\newcommand \cI{{\cal I}}

\newcommand \cK{{\cal K}}

\newcommand \cN{{\cal N}}
\newcommand \cO{{\cal O}}

\newcommand \cS{{\cal S}}

\newcommand \bE{{\mathbb E}}

\newcommand \bP{{\mathbb P}}

\newcommand \bR{{\mathbb R}}

